\documentclass[a4paper,12pt]{amsart}
\usepackage{amsmath,amsthm,amssymb}
\theoremstyle{plain}
\newtheorem{theorem}{Theorem}[section]
\newtheorem{proposition}[theorem]{Proposition}
\newtheorem{lemma}[theorem]{Lemma}

\newtheorem{corollary}[theorem]{Corollary}

\theoremstyle{definition}

\newtheorem{remark}[theorem]{Remark}

\numberwithin{equation}{section}

\newcommand{\field}[1]{\mathbb{#1}}
\newcommand{\C}{\field{C}}

\newcommand{\Q}{\field{Q}}

\newcommand{\Z}{\field{Z}}

\DeclareMathOperator{\Sq}{Sq}

\def\higher{generalized }

\newcommand{\mathscr}{}

\title[Generalized Bott towers]{Topological classification of generalized Bott towers}
\author[S. Choi]{Suyoung Choi}
\address{Department of Mathematical Sciences, KAIST, 335 Gwahangno, Yuseong-gu, Daejeon 305-701, Republic of Korea
} \email{choisy@kaist.ac.kr}
\author[M. Masuda]{Mikiya Masuda}
\address{Department of Mathematics, Osaka City University, Sugimoto, Sumiyoshi-ku, Osaka 558-8585, Japan} \email{masuda@sci.osaka-cu.ac.jp}
\author[D.Y. Suh]{Dong Youp Suh} \address{Department of Mathematical Sciences, KAIST, 335 Gwahangno, Yuseong-gu, Daejeon 305-701, Republic of Korea} \email{dysuh@math.kaist.ac.kr}
\thanks{The first author was partially supported by the second stage of Brain Korea 21 project, KAIST in 2007, the second author was partially supported by Grand-in-Aid for Scientific Research 4102-17540092, and the third author was partially supported by the SRC program of Korea Science and Engineering Foundation R11-2007-035-02002-0.
}
\date{\today}

\begin{document}

\begin{abstract}
If $B$ is a toric manifold and $E$ is a Whitney sum of complex line
bundles over $B$, then the projectivization $P(E)$ of $E$ is again a toric
manifold.
Starting with $B$ as a point and repeating this construction, we obtain
a sequence of complex projective bundles which we call a generalized Bott
tower.
We prove that if the top manifold in the tower has the same cohomology ring
as a product of complex projective spaces, then every fibration in the tower
is trivial so that the top manifold is diffeomorphic to the product of
complex projective spaces.  This gives a supporting evidence to what we call
\emph{cohomological rigidity problem for toric manifolds}
\lq \lq Are toric manifolds diffeomorphic (or homeomorphic)
if their cohomology rings are isomorphic?"
We provide two more results which support the cohomological rigidity
problem.
\end{abstract}

\maketitle

\section{Introduction}\label{Introduction}

A toric variety $X$ of dimension $n$ is a normal complex algebraic
variety with an action of an $n$-dimensional algebraic torus $(\C^*)^n$
having a dense orbit.
A fundamental result in the theory of toric varieties says that there
is a one-to-one correspondence between toric varieties and fans.
It follows that
the classification of toric varieties is equivalent to the
classification of fans up to isomorphism.

Among toric varieties, compact smooth toric varieties, which we
call {\em toric manifolds}, are well studied.  Recently the second
author has shown in \cite{masu07} that toric manifolds as
varieties can be distinguished by their equivariant cohomology. So
we are led to ask how much information \emph{ordinary} cohomology
contains for toric manifolds and posed the following problem in
\cite{ma-su07}. Throughout this paper, an isomorphism of
cohomology rings is as graded rings unless otherwise stated.

\medskip
\noindent
{\bf Cohomological rigidity problem for toric manifolds.}
Are toric manifolds diffeomorphic (or homeomorphic)
if their cohomology rings are isomorphic?

\medskip
If $B$ is a toric manifold and $E$ is a Whitney sum of complex line
bundles over $B$, then the projectivization $P(E)$ of $E$ is again a toric
manifold.
Starting with $B$ as a point and repeating this construction,
say $m$ times, we obtain a sequence of toric manifolds:
\begin{equation} \label{Btower}
B_m\stackrel{\pi_m}\longrightarrow B_{m-1} \stackrel{\pi_{m-1}}\longrightarrow
\dots \stackrel{\pi_2}\longrightarrow B_1 \stackrel{\pi_1}\longrightarrow
B_0=\{\text{a point}\},
\end{equation}
where the fiber of $\pi_i\colon B_i\to B_{i-1}$ for $i=1,\dots,m$ is a
complex projective space $\C P^{n_i}$.  We call the above sequence
a \emph{\higher Bott tower} of height $m$ and omit \lq\lq generalized"
when $n_i=1$ for every $i$ (\cite{gr-ka94}).
We also call $B_k$ in the tower
a $k$-stage \emph{generalized Bott manifold} and omit \lq\lq generalized"
as well when $n_i=1$ for every $i$.

We note that $H^*(B_m)$ is isomorphic to $H^*(\prod_{i=1}^m\C
P^{n_i})$ as a \emph{group} but not necessarily as a \emph{graded
ring}. If every fibration in the tower \eqref{Btower} is trivial,
then $B_m$ is diffeomorphic to $\prod_{i=1}^m\C P^{n_i}$ and
$H^*(B_m)$ is isomorphic to $H^*(\prod_{i=1}^m\C P^{n_i})$ as a graded ring.
The following theorem shows that the converse is
true and generalizes Theorem 5.1 in \cite{MP} treating
Bott towers.
It also gives a supporting evidence to the cohomological rigidity problem
above.

\begin{theorem}\label{main1}
If $H^*(B_m)$ is isomorphic to $H^*(\prod_{i=1}^m\C P^{n_i})$,
then every fibration in the tower \eqref{Btower} is
trivial, in particular, $B_m$ is diffeomorphic to $\prod_{i=1}^m\C
P^{n_i}$.
\end{theorem}

\begin{remark} \label{rem1}
It is shown in \cite{ch-pa-su08} that a toric manifold $X$
whose cohomology ring is isomorphic to that of
a generalized Bott manifold is a generalized Bott manifold.
So we can conclude that if $H^*(X)$ is isomorphic to
$H^*(\prod_{i=1}^m\C P^{n_i})$, then
$X$ is diffeomorphic to $\prod_{i=1}^m\C P^{n_i}$
\end{remark}

2-stage Bott manifolds are famous Hirzebruch
surfaces and their diffeomorphism types can also be distinguished
by their cohomology rings.  The following theorem generalizes this
fact and gives a partial affirmative solution to
the cohomological rigidity problem.

\begin{theorem}[Corollary~\ref{coro:non-singular toric variety}] \label{main2}
2-stage generalized Bott manifolds are diffeomorphic
if and only if their cohomology rings are isomorphic.
\end{theorem}

Actually we obtain a diffeomorphism classification result for
those manifolds (see Theorem~\ref{2stage}) and it would be
interesting to compare it with the variety classification result
in \cite{klei88}.

We also prove the following which gives another partial affirmative solution to
the cohomological rigidity problem.

\begin{theorem}[Theorem~\ref{3stage}] \label{main3}
3-stage Bott manifolds are diffeomorphic
if and only if their cohomology rings are isomorphic.
\end{theorem}

This paper is organized as follows. In Section 2 we recall
well-known facts on projective bundles and discuss their
Pontrjagin classes.  We prepare two lemmas on cohomology of
generalized Bott manifolds in Section 3. Section 4 is devoted to the
proof of Theorem~\ref{main1}. For the proof we need to show that a
Whitney sum of complex line bundles over a product of complex
projective spaces is trivial if its total Chern class is trivial.
This result is of independent interest and
proved in Section 5.
We will discuss 2-stage generalized Bott manifolds in Section 6 and
3-stage Bott manifolds in Section 7.
 In Section 8 which is an appendix, we give a sufficient
condition for an isomorphism of cohomology rings with $\Z/2$ coefficients
to preserve Stiefel-Whitney classes.

\section{Projective bundles}

Let $B$ be a smooth manifold and let $E$ be a complex vector bundle over
$B$.  We denote by $P(E)$ the projectivization of $E$.

\begin{lemma} \label{proje}
Let $B$ and $E$ be as above and let $L$ be a complex line bundle over $B$.
We denote by $E^*$ the complex
vector bundle dual to $E$. Then $P(E\otimes L)$, $P(E)$ and
$P(E^*)$ are isomorphic as fiber bundles over $B$, in particular
they are diffeomorphic.
\end{lemma}

\begin{proof}
For each $x\in B$, we choose a non-zero vector $v_x$ from the fiber of
$L$ over $x$ and define a map $\Psi\colon E\to E\otimes L$ by
$\Psi(u_x):=u_x\otimes v_x$ where $u_x$ is an element of the fiber of $E$
over $x$.
The map $\Psi$ depends on the choice of $v_x$'s but the induced map
from $P(E)$ to $P(E\otimes L)$ does not because $L$ is a line bundle.
It is easy to check that the induced map gives an isomorphism of $P(E)$
and $P(E\otimes L)$ as fiber bundles over $B$.

Choose a hermitian metric $\langle\ ,\ \rangle$ on $E$,
which is anti-$\C$-linear on the first entry and $\C$-linear on
the second entry, and define a map $\Phi\colon E\to E^*$ by
$\Phi(u):=\langle u,\ \rangle$.  This map is not $\C$-linear but
anti-$\C$-linear, so
it induces a map from $P(E)$ to $P(E^*)$, which gives an isomorphism
as fiber bundles.
\end{proof}

Let $y\in H^2(P(E))$ be \emph{minus}\footnote{Our $y$ corresponds to
-$\gamma_1$ in \cite[(2) in p.515]{BH}.} the first Chern class of the
tautological line bundle over $P(E)$ where vectors in a line $\ell$ of $E$
form the fiber over $\ell\in P(E)$.
$H^*(P(E))$ can be viewed as an algebra over $H^*(B)$ via
$\pi^*\colon H^*(B)\to H^*(P(E))$ where $\pi\colon P(E)\to B$ denotes the
projection.  When $H^*(B)$ is finitely generated
and torsion free (this is the case when $B$ is a toric manifold),
$\pi^*$ is injective and
$H^*(P(E))$ as an algebra over $H^*(B)$ is known to be described as
\begin{equation} \label{BoHi}
H^*(P(E))=H^*(B)[y]/\big(\sum_{q=0}^{n}c_q(E)y^{n-q}\big),
\end{equation}
where $n$ denotes the complex fiber dimension of $E$.
If we formally express
\begin{equation} \label{formal}
c(E)=\prod_{i=0}^n(1+u_i),
\end{equation}
then the relation in (\ref{BoHi}) is written as
\begin{equation} \label{uirel}
\sum_{q=0}^{n}c_q(E)y^{n-q}=\prod_{i=0}^n(y+u_i),
\end{equation}
and the total Chern class of the tangent bundle along the fibers $T_f(P(E))$
of $P(E)$ is given by
\begin{equation*} \label{c(TfPE)}
c(T_fP(E))=\prod_{q=0}^{n}(1+y)^{n-q}c_q(E)=\prod_{i=0}^n(1+y+u_i),
\end{equation*}
see \cite[(2) in p.515]{BH}. It follows that the total Pontrjagin class
of $T_f(P(E))$ is given by
\begin{equation} \label{p(TfPE)}
p(T_fP(E))=\prod_{i=0}^n(1+(y+u_i)^2).
\end{equation}

\begin{proposition} \label{Taf}
Let $E'\to B'$ be another complex vector bundle over a smooth manifold $B'$
with the same fiber dimension as $E$.
Suppose that $\varphi\colon H^*(P(E'))\to H^*(P(E))$ is
an isomorphism such that $\varphi(H^*(B'))=H^*(B)$.
Then $\varphi(p(T_fP(E')))=p(T_fP(E))$.
If $\varphi$ satisfies $\varphi(p(B'))=p(B)$ in addition,
then $\varphi(p(P(E')))=p(P(E))$.
\end{proposition}

\begin{proof}
Let $y'$ be an element of $H^2(P(E'))$ defined similarly to $y$.
Since $\varphi$ is an isomorphism and $\varphi(H^*(B'))=H^*(B)$,
we have
\begin{equation} \label{yprime}
\varphi(y')=\epsilon y+w\quad\text{with $\epsilon=\pm 1$ and some
$w\in H^2(B)$}.
\end{equation}
As in (\ref{formal}) we formally express $c(E')=\prod_{i=0}^n(1+u'_i)$.
Then we have the relation (\ref{uirel}) and
the formula (\ref{p(TfPE)}) for $E'\to B'$ with prime.

Since
$\varphi(\prod_{i=0}^n(y'+u'_i))=\prod_{i=0}^n(\epsilon y+w+\varphi(u'_i))$
is zero in $H^*(P(E))$, we have an identity
\begin{equation*} \label{yui}
\prod_{i=0}^n(\epsilon y+w+\varphi(u'_i))=\epsilon^n\prod_{i=0}^n(y+u_i)
\end{equation*}
in a polynomial ring $H^*(B)[y]$ in $y$ with $H^*(B)$ as the coefficient ring.
Replace $y$ with $\sqrt{-1}+y$
and $-\sqrt{-1}+y$ in the identity above and multiply the resulting two
identities at each side.  Then we obtain an identity
\begin{equation} \label{pontui}
\prod_{i=0}^n\big(1+(\epsilon y+w+\varphi(u'_i))^2\big)=\prod_{i=0}^n
\big(1+(y+u_i)^2\big)
\end{equation}
in the ring $H^*(B)[y]$, in particular, in $H^*(P(E))$. It follows from
(\ref{p(TfPE)}), (\ref{yprime}) and (\ref{pontui}) that
\[
\begin{aligned}
\varphi(p(T_fP(E')))&=\varphi\big(\prod_{i=0}^n
(1+(y'+u'_i)^2)\big)\\
&=\prod_{i=0}^n\big(1+(\epsilon y+w+\varphi(u'_i))^2\big)
\\
&=\prod_{i=0}^n\big(1+(y+u_i)^2\big)
\\
&=p(T_fP(E)).
\end{aligned}
\]
This proves the former part of the proposition.

Since the tangent bundle $TP(E)$ of $P(E)$
decomposes into a Whitney sum of $\pi^*(TB)$
and $T_f(P(E))$, we obtain the latter part of the proposition.
\end{proof}

We conclude this section with an observation on Pontrjagin classes
of \higher Bott manifolds in (\ref{Btower}).  Since $\pi_j^*\colon
H^*(B_{j-1})\to H^*(B_j)$ is injective, we regard $H^*(B_{j-1})$
as a subring of $H^*(B_j)$ for each $j$ so that we have a
filtration
\[
H^*(B_m)\supset H^*(B_{m-1})\supset \dots \supset H^*(B_1).
\]

\begin{theorem} \label{filt}
Let \eqref{Btower} be one generalized Bott tower and let
\[
B_m'\to B_{m-1}'\to \dots B_1'\to B_0'=\text{\{a point\}}
\]
be another \higher Bott tower. If $\varphi\colon
H^*(B_m')\to H^*(B_m)$ is an isomorphism which maps
$H^*(B_j')$ onto $H^*(B_j)$ for each $j=1,\dots,m$,
then $\varphi(p(B_j'))=p(B_j)$ for any $j$.
\end{theorem}

\begin{proof}
It follows from the assumption that
the fiber dimensions of $B_j\to B_{j-1}$ and $B_j'\to B_{j-1}'$
must agree for each $j$. If $\varphi(p(B_{j-1}'))=p(B_{j-1})$,
then Proposition~\ref{Taf} implies that $\varphi(p(B_j'))=p(B_j)$.
Therefore, the theorem follows by induction on $j$.
\end{proof}

\section{Cohomology of generalized Bott manifolds} \label{sec:cohom}

Complex vector bundles involved in a generalized Bott tower
(\ref{Btower}) are Whitney sums of complex line bundles.  Since
$P(E\otimes L)$ and $P(E)$ are isomorphic as fiber bundles by
Lemma~\ref{proje}, we may assume that at least one of the complex
line bundles is trivial at each stage of the tower, that is,
$$
B_i=P(\underline{\C}\oplus \xi_i) \quad\text{for $i=1,\dots,m$,}
$$
where $\underline{\C}$ denotes the trivial complex line
bundle and $\xi_i$ a Whitney sum of complex line bundles over $B_{i-1}$.
We set $n_i=\dim \xi_i$.

Let $y_i\in H^2(B_i)$ denote minus the first Chern class of the
tautological line bundle over $B_i=P(\underline{\C}\oplus \xi_i)$.
We may think of $y_i$ as an element of $H^2(B_k)$ whenever $i\le k$.
Then the repeated use of (\ref{BoHi}) shows that the ring
structure of $H^*(B_m)$ can be described as
\begin{equation} \label{BH2}
H^\ast(B_k)=\mathbb Z[y_1,\ldots,y_k]/(f_i(y_1,\ldots,y_i):
i=1,\ldots, k)
\end{equation}
for $k=1,\dots,m$,
where
\begin{equation} \label{yrela}
f_i(y_1,\dots,y_i)=y_i^{n_i+1}+c_1(\xi_i)y_i^{n_i}+\cdots +
c_{n_i}(\xi_i)y_i.
\end{equation}

We prepare two lemmas used later.

\begin{lemma} \label{ym}
The set
$$ \{ by_m+w \in H^2 (B_m) \mid 0\neq b\in\Z,\ w\in H^2(B_{m-1}),
\ (by_m+w)^{n_m +1} = 0\}$$
lies in a one dimensional subspace of $H^2(B_m)$ if it is non-empty.
\end{lemma}

\begin{proof}
We have
\[
\begin{split}
(by_m+&w)^{n_m +1}=(by_m)^{n_m +1} + (n_m +1)(by_m)^{n_m}w +\cdots \\
&=-b^{n_m +1}\sum_{q=1}^{n_m}c_q(\xi_m)y_m^{n_m+1-q}
+(n_m +1)(by_m)^{n_m}w +\cdots
\end{split}
\]
where (\ref{yrela}) is used at the second identity.
If $b \neq 0$ and $(by_m+w)^{n_m+1} = 0$, then
we see $b c_1(\xi_m)=(n_m+1)w$ by looking at the
coefficients of $y_m^{n_m}$ at the
identity above and hence $b$ and $w$ must be proportional, proving the lemma.
\end{proof}

\begin{lemma}\label{power}
Let $x=\sum_{j=1}^mb_jy_j$ be an element of $H^\ast(B_m)$ such that
$b_j\ne 0$ for some $j$. Then $x^{n_j}\ne 0$ in $H^\ast(B_m)$.
\end{lemma}

\begin{proof}
Suppose $x^{n_j}=0$ on the contrary. Then $(\sum_{j=1}^mb_jy_j)^{n_j}$
must be in the ideal generated by the polynomials in (\ref{yrela}) while
a non-zero scalar multiple of $y_j^{n_j}$ appears in
$(\sum_{j=1}^mb_jy_j)^{n_j}$ when we expand it because $b_j\not=0$.
However, it follows from (\ref{yrela}) that if a non-zero scalar multiple
of a power of $y_j$ appears in the ideal, then the exponent must be at least
$n_j+1$, which is a contradiction.
\end{proof}

\section{Cohomologically product generalized Bott manifolds} \label{trivial Bott}

The purpose of this section is to prove Theorem~\ref{main1} in the
Introduction.

We continue to use the notation of the previous section and from
now until this section ends, we assume that $H^*(B_m)$ is
isomorphic to $H^*(\prod_{i=1}^m\C P^{n_i})$. Then, there is
another set of generators $\{x_1,\dots,x_m\}$ in $H^2(B_m)$ such
that
\begin{equation}\label{x1m}
H^\ast(B_m)=\mathbb Z[x_1,\ldots,x_m]/(x_1^{n_1+1},\dots,x_m^{n_m+1}),
\end{equation}
and one has an expression
\begin{equation} \label{yi}
y_i=\sum_{j=1}^m c_{ij}x_j \qquad\textrm{ for }i=1,\ldots, m\textrm{ and }
c_{ij}\in\Z,
\end{equation}
and
\[
x_i=\sum_{j=1}^m d_{ij}y_j \qquad\textrm{ for }i=1,\ldots, m\textrm{ and }
d_{ij}\in\Z,
\]
where both $C=(c_{ij})$ and $D=(d_{ij})$ are unimodular and $C=D^{-1}$.

\begin{lemma} \label{cmm}
By an appropriate change of indices of $x_i$'s with $n_i=n_m$, we may
assume that $c_{mm}=d_{mm}=\pm 1$.
\end{lemma}

\begin{proof}
\emph{Case 1.} The case where all $n_i$'s are same.  In this case
$x_i^{n_m+1}=0$ for any $i$.
Since $x_i=d_{im}y_m+\sum_{j\not=m}d_{ij}y_j$ where the sum lies in
$H^*(B_{m-1})$ by (\ref{BH2}) and $x_i$'s are linearly independent,
it follows from Lemma~\ref{ym} that there is a unique $r$ such that
$d_{rm}\not=0$, and $d_{rm}$ is actually $\pm 1$ because $\det D=\pm 1$.
Since all $n_i$'s are same, we may assume $r=m$ if necessary by changing
the indices of $x_i$'s, so $d_{mm}=\pm 1$ and $d_{im}=0$ for $i\not=m$.
This implies that $c_{mm}=d_{mm}=\pm 1$ and $c_{im}=0$ for $i\not=m$ as well
because $C=D^{-1}$.

\emph{Case 2.} The general case.
Let $S=\{N_1,\ldots, N_k\}$ be the set of all
distinct elements of $n_1,\ldots,n_m$ such that $N_1>\ldots>N_k$.
We can view $\{n_1,\ldots,n_m\}$ as a function
$\mu:\{1,\ldots,m\}\to \mathbb N$ such that $\mu(i)=n_i$. Then $S$
is the image of $\mu$. Let $J_\ell=\mu^{-1}(N_\ell)$ for $\ell=1,\ldots, k$
and let $C_{J_\ell}$ and $D_{J_\ell}$ be the matrices formed from
$c_{ij}$ and $d_{ij}$ with $i,j\in J_\ell$ respectively.

Since $x_i^{n_i+1}=0$, $d_{ij}$ must be $0$ if $n_i< n_j$
by Lemma~\ref{power}.
This shows that $D=(d_{ij})$ is a block upper triangular matrix
$$\left(
\begin{array}{cccc} D_{J_1} & && \ast \\ & D_{J_2} & & \\ & & \ddots& \\ 0
&&&D_{J_k}
\end{array}\right)
$$
if $n_1\ge n_2\ge \dots \ge n_m$, and in general conjugate to the
above by a permutation matrix.

Since $C=D^{-1}$, $C$ is also conjugate to a block upper triangular matrix
$$\left(
\begin{array}{cccc} C_{J_1} & && \ast \\ & C_{J_2} & & \\ & & \ddots& \\ 0
&&&C_{J_k}
\end{array}\right)
$$
by a permutation matrix.
Then a similar argument to Case 1 above can be applied to $C_{J_\ell}$ and
$D_{J_\ell}$ (for $J_\ell$ containing $m$)
instead of $C$ and $D$, and the lemma follows.
\end{proof}

We may further assume that $c_{mm}=d_{mm}=1$ if necessary by taking
$-x_m$ instead of $x_m$, so that we may assume
\begin{equation} \label{ymxm}
y_m=x_m+\sum_{j=1}^{m-1}c_{mj}x_j\quad\text{and}\quad
x_m=y_m+\sum_{j=1}^{m-1}d_{mj}y_j.
\end{equation}

\begin{lemma} \label{prodBk}
If $H^*(B_m)$ is isomorphic to $H^*(\prod_{i=1}^m \mathbb CP^{n_i})$, then
$H^*(B_{m-1})$ is isomorphic to $H^*(\prod_{i=1}^{m-1} \mathbb CP^{n_i})$.
\end{lemma}

\begin{proof}
By (\ref{BH2}),
$H^*(B_{m-1})$ agrees with $H^*(B_m)$ with $y_m=0$ plugged.
On the other hand $H^*(B_m)$ has an expression (\ref{x1m})
by assumption.  It follows from (\ref{ymxm}) that
$H^*(B_{m-1})$ agrees with the right hand side of (\ref{x1m}) with
the relation $x_m+\sum_{j=1}^{m-1}c_{mj}x_j=0$ added.  Therefore,
we can eliminate $x_m$ using the added relation, so that
we obtain a surjective homomorphism
$$\mathbb Z[x_1,\ldots, x_{m-1}]/(x_1^{n_1+1},\dots,x_{m-1}^{n_{m-1}+1})\to
H^\ast(B_{m-1}).$$
But the both sides above are torsion free and have the same rank, so
the homomorphism above is an isomorphism, proving the lemma.
\end{proof}

We need one more result for the proof of Theorem~\ref{main1}.

\begin{theorem}\label{trivial}
A Whitney sum of complex line bundles over a product of complex
projective spaces is trivial if and only if its total Chern class
is trivial.
\end{theorem}

The proof of Theorem~\ref{trivial} is rather long and of
independent interest, so we shall give it in the next section and
complete the proof of Theorem~\ref{main1}.

\begin{proof}[Proof of Theorem~\ref{main1}.]
We prove the theorem by induction on $m$. When $m=1$,
the theorem is obvious.
Assume the theorem is true for $m-1$ case. Suppose $H^*(B_m )
\cong H^* ( \Pi_{j=1}^{m} \C P^{n_j})$. Then
$H^\ast(B_{m-1})\cong H^\ast(\prod_{j=1}^{m-1}\mathbb CP^{n_j})$
by Lemma~\ref{prodBk}. By the induction hypothesis,
the \higher Bott tower (\ref{Btower}) is trivial up to $B_{m-1}$,
in particular,
$B_{m-1}$ is diffeomorphic to $\prod_{j=1}^{m-1}\mathbb CP^{n_j}$.

Remember that $B_m=P(\underline{\C}\oplus \xi_m)$ and express
$\xi_m=\bigoplus_{i=1}^{n_m}\eta_i$ where $\eta_i$ is a complex line
bundle over $B_{m-1}$.
Let $\gamma_j$ be the
complex line bundle over $B_{m-1}$ whose first Chern class is $y_j$.
One can write
\begin{equation} \label{c1eta}
c_1(\eta_i)=\sum_{j=1}^{m-1}a_{ij}y_j \quad\text{with $a_{ij}\in\mathbb Z$}.
\end{equation}
Then $\eta_i=\bigotimes_{j=1}^{m-1}\gamma_j^{a_{ij}}$.

The fibration $B_m=P(\underline{\mathbb C}\oplus\xi_m)\to B_{m-1}$ is
isomorphic to
a fibration $P((\underline{\mathbb C}\oplus \xi_m)\otimes L) \to B_{m-1}$
for any complex line bundle $L$ over $B_{m-1}$.
Therefore, it suffices to find a complex line bundle $L$ such that the total
Chern class of $(\underline{\mathbb C}\oplus \xi_m)\otimes L$ is trivial
because
the triviality of the bundle follows from the triviality of the Chern class
by Theorem~\ref{trivial}.

We take
$L=\bigotimes_{j=1}^{m-1}\gamma_j^{-d_{mj}}$ with $d_{mj}$ in
(\ref{ymxm}).  Then
\begin{equation} \label{chern}
c((\underline{\mathbb C}\oplus \xi_m)\otimes L)
=\prod_{i=0}^{n_m} \left(
            1+ \sum_{j=1}^{m-1} (a_{ij}-d_{mj})y_j \right)\quad\text{in
            $H^*(B_{m-1})$,}
\end{equation}
where $a_{0j}=0$.
On the other hand, it follows from (\ref{BH2}) and (\ref{c1eta}) that
\begin{eqnarray}
f_m(y_1,\dots,y_m)&=&y_m^{n_m+1}+c_1(\xi_m)y_m^{n_m}+\cdots +
c_{n_m}(\xi_m)y_m \nonumber\\
& = & y_m\prod_{i=1}^{n_m}(y_m+c_1(\eta_i))\nonumber\\
& = & \prod_{i=0}^{n_m}(y_m+\sum_{j=1}^{m-1}a_{ij}y_j)\label{eq:product form of f_m}\\
& = & 0\nonumber
\end{eqnarray}
in $H^\ast (B_m)$. We plug $y_m = x_m-\sum_{j=1}^{m-1}d_{mj}y_j$ from
(\ref{ymxm}) into (\ref{eq:product form of f_m}) to get
\begin{equation} \label{0xm}
0
    = \prod_{i=0}^{n_m}(x_m + \sum_{j=1}^{m-1}(a_{ij}-d_{mj})y_j)
\end{equation}
in $H^*(B_m)$. Here we note that
\[
H^*(B_m)=\Z[y_1,\dots,y_{m-1},x_m]/(f_1(y_1),\dots,f_{m-1}(y_1,\dots,y_{m-1}),
x_m^{n_m+1})
\]
because a natural homomorphism from the right hand side above to
$H^*(B_m)$ is surjective by (\ref{BH2}) and (\ref{ymxm}), and hence isomorphic
since both are torsion free and have the same rank.  Therefore,
when we expand the right hand side of (\ref{0xm}), the coefficient of
$x_m^k$ must be zero
for any $k=1,\dots,n_m$.  This implies that the right hand side of
(\ref{chern}) is equal to $1$, proving the theorem.
\end{proof}

Combining Theorem~\ref{main1} with Theorem 8.1 in \cite{ch-ma-su07}, we obtain
the following corollary which generalizes Theorem 5.1 in \cite{MP}
treating the case where $n_i=1$ for any $i$.

\begin{corollary}
If the cohomology ring of a quasitoric manifold over a product of
simplices is isomorphic to
that of $\prod_{i=1}^m\C P^{n_i}$, then it is homeomorphic to
$\prod_{i=1}^m\C P^{n_i}$.
\end{corollary}

\begin{remark}
Similarly to Remark~\ref{rem1} the assumption
\lq\lq over a product of simplices" in the corollary above can be dropped
by a result in \cite{ch-pa-su08}.
\end{remark}

\section{Proof of Theorem~\ref{trivial}}

This section is devoted to the proof of Theorem~\ref{trivial}.
We recall a general fact.
A more refined result can be found in \cite{Pe}.

\begin{lemma} \label{stable}
Let $X$ be a finite CW-complex such that $H^{odd}(X)=0$ and
$H^*(X)$ has no torsion.
Then complex $n$-dimensional vector bundles over $X$ with
$2n\ge \dim X$ are isomorphic if and only if
their total Chern classes are same.
\end{lemma}

\begin{proof}
The assumption on $H^*(X)$ implies that $K(X)$ is torsion free, so Chern
character gives a monomorphism from $K(X)$ to $H^{*}(X;\Q)$.
On the other hand, if $\dim X\le 2n$, then the homotopy set $[X,BU(n)]$,
where $BU(n)$ denotes the classifying space of a unitary group $U(n)$,
agrees with $K(X)$.  This implies the lemma.
\end{proof}

Let $B=\prod_{j=1}^k \C P^{n_j}$ be a product of complex
projective spaces and let $E=\bigoplus_{i=1}^{n}\eta_i$ be a Whitney sum
of complex line bundles $\eta_i$ over $B$. Suppose that $c(E)=1$.
Then since $H^{odd}(B)=0$ and $H^*(B)$ has no torsion, $E$ is trivial by
Lemma~\ref{stable} when $n\ge \sum_{j=1}^{k}n_j$. So we assume
that $n<\sum_{j=1}^{k}n_j$ in the following.

By assumption
$$H^\ast(B)=\mathbb Z[x_1,\ldots,x_{k}]/(x_1^{n_1+1},\dots,x_k^{n_k+1}),$$
where we can take $x_j$ as the first Chern class of the pullback $\gamma_j$
of the tautological line bundle over $\C P^{n_j}$ via the projection
$\prod_{j=1}^{k}\mathbb CP^{n_j}\to \mathbb CP^{n_j}$.
Then we may assume that $\eta_i=\bigotimes_{j=1}^{k}\gamma_j^{a_{ij}}$ with
$a_{ij}\in\Z$ and
\begin{equation} \label{cE}
1=c(E)=\prod_{i=1}^n(1+\sum_{j=1}^{k}a_{ij}x_j).
\end{equation}
It follows that
$$0=c_1(E)=\sum_{i=1}^n(\sum_{j=1}^{k}a_{ij}x_j)=\sum_{j=1}^{k}
\left(\sum_{i=1}^n a_{ij}\right)x_j.
$$
Since $x_j$'s are linearly independent, the identity above
implies that
\begin{equation}\label{eq:sum a_{ij}=0}
\sum_{i=1}^n a_{ij}=0\qquad \textrm{for each  \ }j=1,\ldots, k.
\end{equation}
Moreover it follows from (\ref{cE}) that
\begin{equation}\label{eq:second chern class=0}
0=c_2(E)=\sum_{i'>i=1}^n\left[\left(\sum_{j=1}^{k}a_{ij}x_j\right)
\left(\sum_{j=1}^{k}a_{i'j}x_j\right)\right].
\end{equation}

We need to consider two cases.

\noindent \underline{\textbf{Case I} \quad $n_j\ge 2$ for some
$j=1,\ldots,k$.}

Since $x_j^2\ne 0$ in $H^\ast(B)$ in this case, the coefficient of
$x_j^2$-term in (\ref{eq:second chern class=0}) must vanish. Thus
$\sum_{i'>i=1}^n a_{ij}a_{i'j}=0$. Therefore from (\ref{eq:sum
a_{ij}=0}) we have
$$0=\left(\sum_{i=1}^n a_{ij}\right)^2=\sum_{i=1}^na_{ij}^2 +
2\sum_{i'>i=1}^n a_{ij}a_{i'j}.$$ Hence $\sum_{i=1}^na_{ij}^2=0$,
which implies that
$$a_{1j}=\cdots=a_{nj}=0.$$

\noindent \underline{\textbf{Case II} \quad $n_j=1$ for all
$j=1,\ldots,k$.}

In this case $n<k$ as $n<\sum_{j=1}^{k}n_j$. Set $\mathbf
v_j=(a_{1j},\ldots, a_{nj})\in \mathbb Z^n$ for $j=1,\ldots, k$.
We claim that $\mathbf v_j=\mathbf 0$ for some $j=1,\ldots, k$.
Since $x_jx_{j'}\ne 0$ in $H^\ast(B)$ for $j\ne j'$, the coefficient
of $x_jx_{j'}$-term in (\ref{eq:second chern class=0}) must vanish.
Namely,
\begin{eqnarray*}
0 & = & \sum_{i=1}^n a_{ij}\left(\sum_{i'=1}^n
a_{i'j'}-a_{ij'}\right)\\
& = & \left(\sum_{i=1}^n
a_{ij}\right)\left(\sum_{i'=1}^na_{i'j'}\right)-\sum_{i=1}^na_{ij}a_{ij'}.
\end{eqnarray*}
By (\ref{eq:sum a_{ij}=0}) we have
$$\sum_{i=1}^na_{ij}a_{ij'}=0\quad \textrm{for all }1\le j, j'\le k,
\quad j\ne j'.$$ This means that $k$ vectors $\mathbf v_1,
\ldots, \mathbf v_{k}$ in $\mathbb Z^n\subset \mathbb R^n$ are
mutually orthogonal. But since $k>n$, $\mathbf v_j=\mathbf 0$ for some
$j=1,\ldots, k$.

We have shown that in either cases there exists some $j$ such that
$(a_{1j},\ldots, a_{nj})=\mathbf 0$. For simplicity,
assume $j=k$. Then $\eta_i$ is of
the form $\bigotimes_{j=1}^{k-1}\gamma_j^{a_{ij}}$. Let
$\overline{\gamma}_j$ be the pull-back bundle of the tautological line
bundle of $\mathbb CP^{n_j}$ via the projection
$\prod_{j=1}^{k-1}\mathbb CP^{n_j}\to \mathbb CP^{n_j}$. Then
$E=\left(\bigoplus_{i=1}^n\bigotimes_{j=1}^{k-1}
\overline{\gamma}_j^{a_{ij}}\right)\times
\mathbb CP^{n_{k}}$. Hence the problem reduces to the bundle on
$\prod_{j=1}^{k-1}\mathbb CP^{n_j}$.

The argument above shows that the proof of the theorem reduces to the case
$k=1$, so the theorem follows from the following lemma.

\begin{lemma} \label{rigid}
Let $E$ and $E'$ be Whitney sums of complex lines bundles over $\C
P^n$ of the same dimension. If $c(E)=c(E')$, then $E$ and $E'$ are
isomorphic.
\end{lemma}

\begin{proof}
Let $\gamma^u$ denote a complex line bundle over $\C P^n$ whose
first Chern class is $u\in H^2(\C P^n)$.  Then
$E=\oplus_{i=0}^m\gamma^{u_i}$ and
$E'=\oplus_{i=0}^m\gamma^{u_i'}$ with $u_i,u_i'\in H^2(\C P^n)$.
In case $m\ge n$, the lemma follows from Lemma~\ref{stable}. In
case $m<n$, $c(E)=c(E')$ implies that
$\{u_0,\dots,u_m\}=\{u_0',\dots,u_m'\}$ and hence $E$ and $E'$ are
isomorphic.
\end{proof}


\section{2-stage generalized Bott manifolds}

A 2-stage Bott manifold is a Hirzebruch surface
$H_a$ and it is well-known that $H_a$ and
$H_b$ are isomorphic as varieties if and only if $|a|=|b|$ and
diffeomorphic if and only if $a\equiv b\pmod 2$.

2-stage \higher Bott manifolds can be thought of
as a higher dimensional generalization of Hirzebruch surfaces and
their classification as varieties is completed in \cite{klei88}.
In this section we complete the diffeomorphism classification of
those manifolds.

Let $B_1=\C P^{n_1}$ and
\[
B_2=P(\bigoplus_{i=0}^{n_2}\gamma^{u_i}),
\]
where $u_0=0$ and
$\gamma^{u_i}$ denotes the complex line bundle over $B_1$ whose
first Chern class is $u_i\in H^2(B_1)$ as before.
Similarly let
\[
B_2'=P(\bigoplus_{i=0}^{n_2}\gamma^{u_i'})
\]
be another 2-stage \higher Bott manifold with
$B_1=\C P^{n_1}$ as 1-stage, where $u_0'=0$.

\begin{theorem}\label{2stage}
Let $B_2$ and $B_2'$ be as above. Then the following are equivalent.
\begin{enumerate}
\item There exist $\epsilon=\pm 1$ and $w\in H^2(B_1)$ such that
$$\prod_{i=0}^{n_2}(1+\epsilon (u_i'+w))=\prod_{i=0}^{n_2}(1+u_i)
    \quad\text{in $H^*(B_1)$}.$$
\item $B_2$ and $B_2'$ are diffeomorphic.
\item $H^*(B_2)$ and $H^*(B_2')$ are isomorphic.
\end{enumerate}
\end{theorem}

\begin{proof}
Condition (1) means that
$(\bigoplus_{i=0}^{n_2}\gamma^{u_i'})\otimes \gamma^w$ or
its dual has
the same total Chern class as $\bigoplus_{i=0}^{n_2}\gamma^{u_i}$,
so that they are isomorphic as vector bundles by Lemma~\ref{rigid}.
This together with Lemma~\ref{proje} implies (2).
The implication (2)$\Rightarrow$(3) is obvious, so it suffices to prove
the implication (3)$\Rightarrow$(1).

Suppose $H^\ast(B_2)$ and $H^\ast(B_2')$ are isomorphic.
Then there is an isomorphism
\[
\begin{split}
\varphi\colon H^*(B_2')=\Z[x,y']/&(x^{n_1+1},
\prod_{i=0}^{n_2}(y'+u_i'))\\
\to &H^*(B_2)=\Z[x,y]/(x^{n_1+1}, \prod_{i=0}^{n_2}(y+u_i)).
\end{split}
\]
Express
\begin{equation} \label{eqn:XY}
\varphi(x)=px+qy \quad\text{and}\quad \varphi(y')=rx+sy
\end{equation}
with $p,q,r,s\in \mathbb Z$.
Since $\varphi$ is an isomorphism, we have
\begin{equation}\label{determinant for X and Y}
ps-qr=\pm1.
\end{equation}

We distinguish three cases.

{\it Case 1.} The case where $n_1\ge 2$ and $n_2=1$. We write $u_1=ax$
and $u_1'=a'x$.
Since $\varphi(y'(y'+a'x))=0$ and $y(y+ax)=0$ in $H^\ast(B_2)$, we have
\begin{eqnarray*}
0&=&(rx+sy)((rx+sy)+a'(px+qy))\\
&=&r(r+a'p)x^2 + \big(r(s+a'q) + s(r+a'p)\big)x y + s(s+a'q)y^2\\
&=&r(r+a'p)x^2+ \big(r(s+a'q) + s(r+a'p)-s(s+a'q)a\big)xy.
\end{eqnarray*}
Therefore,
\begin{equation}\label{relation of a and a'}
r(s+a'q) + s (r+a'p)=s(s+a'q)a
\end{equation}
and moreover since $n_1\ge 2$, we have $r(r+a'p)=0$ and hence
$r=0$ or $r=-a'p$.

If $r=0$, then $p=\pm1$ and $s=\pm1$ from (\ref{determinant for X and Y})
and hence $\pm a'=(s+a'q)a$ from (\ref{relation of a and a'}),
which implies that $a|a'$.
If $r=-a'p$, then from (\ref{determinant for X and Y}) we have
$\pm1=ps-qr=ps+a'pq=p(s+a'q)$. Thus $p=\pm1$ and $s+a'q=\pm1$.
From (\ref{relation of a and a'}) we have
$\pm a'=sa$ and hence $a | a'$.
In any case we have shown that $a'$ is divisible by $a$.
By the symmetry, $a$ is divisible by $a'$.
Thus $a=\pm a'$, and hence the identity in (1) is satisfied with $w=0$.

{\it Case 2.} The case where $n_1=n_2=1$. We write $u_1=ax$ and $u_1'=a'x$
as in Case 1 above. The identity in (1) is equivalent to
$$a\equiv a' \mod 2.$$
In the following all congruence relations are taken modulo 2 unless stated
otherwise.  It follows from (\ref{relation of a and a'}) and
(\ref{determinant for X and Y}) that
\begin{equation} \label{eqn:1}
a'\equiv s(s+a'q)a.
\end{equation}
On the other hand, since $x^2=0$, the identity $\varphi(x)^2=0$ implies that
\[
0=(px+qy)^2\equiv q^2y^2\equiv q^2axy,
\]
so that
\begin{equation} \label{eqn:2}
q^2a\equiv 0.
\end{equation}
If $a\equiv 0$, then so is $a'$ from (\ref{eqn:1}).  If $a\equiv 1$,
then $q\equiv 0$ from (\ref{eqn:2}), so that $a'\equiv s^2a$
from (\ref{eqn:1}) and $ps\equiv 1$ (and hence $s\equiv 1$)
from (\ref{determinant for X and Y}).
Therefore, $a\equiv a'$ in any case.

{\it Case 3.} The case where $n_2 \ge 2$. If $u_i=u_i'=0$ for all
$i$'s, then the identity in (1) holds with $w=0$, so we may assume
either $u_i$ or $u'_i$ is non-zero for some $i$ and moreover
$u_i\not=0$ for some $i$ without loss of generality. Then, since
$0=\varphi(x^{n_1+1})=(px+qy)^{n_1+1}$, $q=0$ by
Lemma~\ref{lemma:technical} below. This means that $\varphi$
preserves the subring $H^*(B_1)$, so that we have
\[
\varphi(y')=\epsilon y+w \quad\text{for some $w\in H^2(B_{1})$},
\]
where $\epsilon =\pm 1$. Therefore
$\varphi(\prod_{i=0}^{n_2}(y'+u'_i))= \prod_{i=0}^{n_2}(\epsilon
y+w+\varphi(u'_i)).$ Since this element vanishes in $H^*(B_2)$ and
is a polynomial of degree $n_2+1$ in $y$, we have an identity
\[
\prod_{i=0}^{n_2}(y+\epsilon
(w+\varphi(u'_i)))=\prod_{i=0}^{n_2}(y+u_i)
\]
as polynomials in $y$.  Then, plugging $y=1$, we obtain the
identity in (1) in the theorem.
\end{proof}

Here is the lemma used above. We shall use the same notation as above.

\begin{lemma}\label{lemma:technical}
Assume that $n_2\ge 2$, $u_0=0$ and $u_i\ne 0$ for some $1\le i\le n_2$.
If $(\alpha x+\beta y)^{n_1+1}=0$
in the ring $\mathbb Z[x,y]/(x^{n_1+1}, \prod_{i=0}^{n_2}(y+u_i))$
for some integers $\alpha, \beta$, then $\beta=0$.
\end{lemma}

\begin{proof}
Since $(\alpha x+\beta y)^{n_1+1}=0$ in the ring $\mathbb Z[x,y]/(x^{n_1+1},
\prod_{i=0}^{n_2}(y+u_i))$, there are a homogeneous polynomial $g(x,y)$
in $x,y$ of total degree $n_1-n_2$ and an integer $c$ such that
\begin{equation}\label{equation:polynomial on $x$ and $y$}
(\alpha x+\beta y)^{n_1+1}-cx^{n_1+1}=g(x,y)\prod_{i=0}^{n_2}(y+u_i)
\end{equation}
as polynomials in $x$ and $y$.  In fact, $c=\alpha^{n_1+1}$ as $u_0=0$.

Suppose $g(x,y)\ne 0$ (so that $n_1\ge n_2$).
When we split the left hand side into a product of linear polynomials
in $x$ and $y$, it has at most two linear polynomials over $\Z$ as factors
while the right hand side has at least three linear polynomials over $\Z$ as
$n_2\ge 2$ by assumption.  This is a contradiction.
Therefore $g(x,y)=0$. But then $\beta$ must be zero, proving the lemma.
\end{proof}

\begin{corollary}\label{coro:non-singular toric variety}
2-stage \higher Bott manifolds are
diffeomorphic if and only if their cohomology rings
are isomorphic.
\end{corollary}

\begin{proof}
Let $B_2\to B_1=\C P^{n_1}$ be a \higher Bott tower of height 2 where
the fiber is $\C P^{n_2}$ and let $B_2'\to B_1'=\C P^{n_1'}$ be
another \higher Bott tower of height 2 where the fiber is $\C P^{n_2'}$.
Suppose that $H^*(B_2)$ is isomorphic to $H^*(B_2')$.  Then
$\{n_1,n_2\}=\{n_1',n_2'\}$ which we can see from their Betti
numbers. If $n_i=n_i'$ for $i=1,2$, then the corollary follows
from Theorem~\ref{2stage}.  Therefore, we may assume that
$n_1=n_2'$, $n_2=n_1'$ and they are different. If both $B_2$ and
$B_2'$ are cohomologically product, then they are diffeomorphic to
$\C P^{n_1}\times \C P^{n_2}$ by Theorem~\ref{main1}.

In the sequel it suffices to prove that
$H^*(B_2)$ and $H^*(B_2')$ are not isomorphic when they are not
cohomologically product and $n_1=n_2'\not=n_2=n_1'$.
We may assume $n_1>n_2$ without loss of generality.
Since $B_2'$ is a
$\mathbb CP^{n_1}$-bundle over $\mathbb CP^{n_2}$, there is a non-zero element
in $H^2(B_2')$ whose $n_1$-th power vanishes, in fact,
a non-zero element in $H^2(B_2')$ coming from the base space
$\mathbb CP^{n_2}$ is such an element because $n_1>n_2$.
On the other hand, it is not difficult to see that there is no such
a non-zero element in $H^2(B_2)$ since
\[
H^*(B_2)=\Z[x,y]/(x^{n_1+1},\prod_{i=0}^{n_2}(y+u_i)),
\]
where $u_0=0$ and $u_i\ne 0$ for some $1\le i\le n_2$. (It also follows from
Lemma~\ref{lemma:technical} when $n_2 \ge 2$.)
\end{proof}

\section{3-stage Bott manifolds}

This section is devoted to the proof of Theorem~\ref{main3} in the Introduction, that is

\begin{theorem} \label{3stage}
3-stage Bott manifolds are diffeomorphic
if and only if their cohomology rings are isomorphic.
\end{theorem}

Remember a Bott tower of height $3$
\[
B_3\stackrel{\pi_3}\longrightarrow
B_2\stackrel{\pi_2}\longrightarrow
B_1\stackrel{\pi_1}\longrightarrow B_0=\{\text{a point}\}
\]
where $B_i=P(\underline{\C}\oplus\xi_i)$ for $i=1,2,3$ and
$\xi_i$ is a complex line bundle over $B_{i-1}$.
Let $\gamma_i$ be the dual of the tautological line bundle over $B_{i}$
and let $y_i$ be the first Chern class of $\gamma_i$.  Then
\[
\xi_1=\underline{\C},\quad
\xi_2=\gamma_1^a,\quad \xi_3=\pi_2^*(\gamma_1)^b\otimes\gamma_2^c
\]
with integers $a,b,c$, and it follows from (\ref{BoHi}) that
\begin{equation} \label{HB3}
H^*(B_3)=\Z[y_1,y_2,y_3]/\big(y_1^2, y_2(ay_1+y_2), y_3(by_1+cy_2+y_3)\big),
\end{equation}
where $y_1$ and $y_2$ are regarded as elements of $H^*(B_3)$ as before.

We note that $H^*(B_2;\Q)$, that is a subring of $H^*(B_3;\Q)$ with $y_3=0$
in \eqref{HB3}, is isomorphic to $H^*((\C P^1)^2;\Q)$ because $y_1^2=0$ and
$(ay_1+2y_2)^2=0$.  However $H^*(B_3;\Q)$ is not necessarily isomorphic to
$H^*((\C P^1)^3;\Q)$ as is shown in the following lemma.

\begin{lemma} \label{B3CP13}
The following are equivalent.
\begin{enumerate}
\item $H^*(B_3;\Q)\cong H^*((\C P^1)^3;\Q)$.
\item $(\sum_{i=1}^3a_iy_i)^2=0$ in $H^*(B_3;\Q)$ for some integers (or rational numbers)
$a_i$ with $a_3\not=0$.
\item $c(2b-ac)=0$.
\end{enumerate}
\end{lemma}

\begin{proof}
(1) $\Leftrightarrow$ (2). This equivalence follows from the observation made
in the paragraph just before the lemma.

(2) $\Rightarrow$ (3).
Since $y_2^2=-ay_1y_2$ and $y_3^2=-by_1y_3-cy_2y_3$ by (\ref{HB3}), we have
\[
\begin{split}
(\sum_{i=1}^3&a_iy_i)^2=\sum_{i=1}^3a_i^2y_i^2+2\sum_{i<j}a_ia_jy_iy_j\\
&=(2a_1a_2-a_2^2a)y_1y_2+(2a_2a_3-a_3^2c)y_2y_3+(2a_3a_1-a_3^2b)y_3y_1,
\end{split}
\]
so that
\[
2a_1a_2=a_2^2a, \quad 2a_2a_3=a_3^2c, \quad 2a_3a_1=a_3^2b.
\]
An elementary computation shows that these imply $c(2b-ac)=0$.

(3) $\Rightarrow$ (2).
If $c(2b-ac)=0$, then an elementary computation shows that
$(by_1+cy_2+2y_3)^2=0$.
\end{proof}

\begin{lemma} \label{p1B3}
The first Pontrjagin class of $B_3$ is given by
\[
p_1(B_3)=c(2b-ac)y_1y_2.
\]
Therefore, $p_1(B_3)=0$ if and only if $H^*(B_3;\Q)\cong
H^*((\C P^1)^3;\Q)$.
\end{lemma}

\begin{proof}
Since $p(B_1)=1$, it follows from (\ref{p(TfPE)}) that
\[
p(B_2)=(1+y_2^2)(1+(y_2+ay_1)^2)
\]
which is equal to $1$ because $y_2^2=-ay_1y_2$.  Therefore,
it follows from (\ref{p(TfPE)}) again that
\[
\begin{split}
p(B_3)&=(1+y_3^2)(1+(y_3+by_1+cy_2)^2)\\
&=1+y_3^2+(y_3+by_1+cy_2)^2\\
&=1+c(2b-ac)y_1y_2
\end{split}
\]
where we used $y_3(y_3+by_1+cy_2)=0$ and $y_2^2=-ay_1y_2$. This proves the
former part of the lemma. The latter part follows from the former part and
Lemma~\ref{B3CP13}.
\end{proof}

We shall complete the proof of Theorem~\ref{3stage}.
Let $B_3'\to B_2'\to B_1'$ be another Bott tower of height $3$ and we denote
the elements corresponding to $a,b,c$ and $y_i$ by $a',b',c'$ and $y_i'$.
The results in \cite{Wa} and \cite{Ju}
tell us that if there is an isomorphism $\varphi\colon H^*(B_3')
\to H^*(B_3)$ such that $\varphi(p_1(B_3'))=p_1(B_3)$ and $\varphi(w_2(B_3'))
=w_2(B_3)$, then $B_3$ and $B_3'$ are diffeomorphic, where $w_2$ denotes the
second Stiefel-Whitney class.

Suppose $H^*(B_3')\cong H^*(B_3)$ and let $\varphi\colon H^*(B_3')
\to H^*(B_3)$ be an isomorphism.  Since $\varphi(w_2(B_3'))=w_2(B_3)$ by
Lemma~\ref{SW} in the appendix,
it suffices to check $\varphi(p_1(B_3'))=p_1(B_3)$.

If $H^*(B_3';\Q)\cong H^*(B_3;\Q)$ is isomorphic to $H^*((\C P^1)^3;\Q)$,
then $p_1(B_3')=p_1(B_3)=0$ by Lemma~\ref{p1B3}, in particulart
$\varphi(p_1(B_3'))=p_1(B_3)$. Suppose that
$H^*(B_3';\Q)\cong H^*(B_3;\Q)$ is not
isomorphic to $H^*((\C P^1)^3;\Q)$.  Then Lemma~\ref{B3CP13} says that
there is no element
$\sum_{i=1}^3a_iy_i$ for rational numbers $a_i$'s with $a_3\not=0$ such that
$(\sum_{i=1}^3a_iy_i)^2=0$.  On the other hand,
$y_1^2=(a/2y_1+y_2)^2=0$ and $y_1$ and $a/2y_1+y_2$ generate the subring
$H^*(B_2;\Q)$.  The same holds for $B_3'$. It follows that the images of
$y_1'$ and $a'/2y_1'+y_2'$ by $\varphi$ generate the subring $H^*(B_2;\Q)$,
and hence $\varphi(H^*(B_2'))\subset H^*(B_2)$.  Therefore
Proposition~\ref{Taf} can be applied and we conclude that
$\varphi(p_1(B_3'))=p_1(B_3)$ because $p(B_2')=p(B_2)=1$.

\section{Appendix}

In this appendix, we prove a general fact used in the previous section on
Stiefel-Whitney classes.
In the following, cohomology will be taken with $\Z/2$ coefficients
unless otherwise stated.
Let $M$ be a connected closed manifold of dimension $n$ and let
\[
\Sq(x)=x+\Sq^1(x)+\Sq^2(x)+\dots+\Sq^n(x)\quad \text{for $x\in H^*(M)$}
\]
denote the total squaring operation, where $\Sq^k\colon H^q(M)\to
H^{q+k}(M)$ is an additive homomorphism.
The $k$-th Wu class $v_k(M)\in H^k(M)$ of $M$ is characterized by
\begin{equation} \label{Wu}
v_k(M)\cup x=\Sq^k(x) \quad \text{for any $x\in H^{n-k}(M)$}
\end{equation}
and
\begin{equation} \label{sqw}
\Sq(v(M))=w(M)
\end{equation}
where $v(M)=1+v_1(M)+v_2(M)+\dots+v_n(M)$ and $w(M)$ denotes the total
Stiefel-Whitney class of $M$, see \cite[p.132]{mi-st74}.


\begin{lemma} \label{SW}
Suppose that $H^*(M)$ is generated by $H^r(M)$ for some $r$ as a ring and
let $M'$ be another connected closed manifold of dimension $n$ such that
$H^*(M')$ is isomorphic to $H^*(M)$ as a ring.
Then $\phi(w(M'))=w(M)$ for any ring isomorphism $\phi\colon H^*(M')\to H^*(M)$.
\end{lemma}

\begin{proof}
We first show that $\phi$ commutes with $\Sq$.  Since $H^*(M')\cong H^*(M)$
are generated by elements of degree $r$, $\Sq(y)=y+y^2$ and
$\Sq(\phi(y))=\phi(y)+\phi(y)^2$ for $y\in H^r(M')$; so
$\phi(\Sq(y))=\Sq(\phi(y))$.  This implies $\phi(\Sq(y))=\Sq(\phi(y))$
for any $y\in H^*(M')$ because both $\phi$ and
$\Sq$ are ring homomorphisms and $H^*(M')\cong H^*(M)$ are generated by
elements of degree $r$ as rings.

It follows from (\ref{Wu}) and the commutativity of $\phi$ and $\Sq$ that
\[
\phi(v_k(M'))\cup \phi(y)=\phi(\Sq^k(y))=\Sq^k(\phi(y))
\]
for any $y\in H^{n-k}(M')$. Since $\phi(H^{n-k}(M'))=H^{n-k}(M)$,
the above identity together with (\ref{Wu}) implies
\begin{equation} \label{phiv}
\phi(v_k(M'))=v_k(M).
\end{equation}
It follows from (\ref{sqw}), (\ref{phiv}) and the commutativity of $\phi$
and $\Sq$ that
\[
\phi(w(M'))=\phi(\Sq(v(M')))=\Sq(\phi(v(M')))=\Sq(v(M))=w(M),
\]
proving the lemma.
\end{proof}

\end{document}